\documentclass[12pt]{amsart}
\usepackage{a4wide,times,graphics,amssymb,psfrag}

\title[Separation of singularities]
{Uniform estimates in the Poincar\'e-Aronszajn theorem on the 
separation of singularities of analytic functions}
\author{V.P. Havin}
\address{St. Petersburg State University. Department of Mathematics and
Mechanics, 28 Universitetski pr., St. Petersburg, 198504, Russia}
\email{havin@VH1621.spb.edu}
\author{A.H. Nersessian}
\address{Burnside Hall, Department of Mathematics and Statistics McGill
University, Montreal H3A 2K6, Canada}
\author{J. Ortega-Cerd\`a}
\address{Dept.\ Matem\`atica Aplicada i An\`alisi, Universitat  de Barcelona,
Gran Via 585, 08071 Bar\-ce\-lo\-na, Spain}
\email{jortega@ub.edu}
\thanks{The first author is partially supported by RFBR grant \#03-01-00377
and by the grant for Leading Scientific Schools \#NSH-2266-2003.1. The last
author is supported by DGICYT grant BFM2002-04072-C02-01 and the CIRIT grant
2001-SGR00172}

\keywords{} 
\subjclass{}
\date{\today}

\newcommand{\D}{\mathbb D}
\newcommand{\C}{\mathbb C}
\newcommand{\R}{\mathbb R}

\newcommand{\K}{\mathcal K}

\newcommand{\dist}{\operatorname{dist}}
\newcommand{\Clos}{\operatorname{Clos}}
\newcommand{\Arg}{\operatorname{arg}}
\newcommand{\dbar}{\bar\partial}

\newtheorem{theorem}{Theorem}

\newtheorem{lemma}{Lemma}
\newtheorem*{definition}{Definition}
\theoremstyle{definition}

\begin{document}
\begin{abstract}
We study the possibility of splitting any bounded analytic
function $f$ with singularities in a closed set $E\cup F$ as a sum of
two bounded analytic functions with singularities in $E$ and $F$ respectively.
We obtain some results under geometric restrictions on
the sets $E$ and $F$ and we provide some examples showing the sharpness
of the positive results.
\end{abstract}
\maketitle
\section*{Introduction}
Let $O\subset \C$ be an open set and $S_1,S_2$ be two relatively closed subsets
in
$O$, $S=S_1\cup S_2$. The following result is due to Aronszajn
\cite{Aronszajn35}
\begin{theorem}\label{Aron}
Any function $f$ analytic in $O\setminus S$ coincides with
$(f_1+f_2)|_{(O\setminus S)}$ where $f_j$ are analytic in $O\setminus S_j$,
$j=1,2$.
\end{theorem}
In modern textbooks this fact is treated (if at all) as a trivial example
illustrating general sheaf theoretic and $\bar\partial$ approaches and related
to the first Cousin problem; \cite{Aronszajn35} is never quoted (see e.g.
\cite[p.\ 225]{BerGay91}, \cite{Hormander90a})  not to mention its famous
predecessors Poincar\'e \cite{Poincare92} and Fr\'echet \cite{Frechet30} 
although their approaches remain interesting even now. The first version of
Theorem~\ref{Aron} with an ingenious proof appeared in 1892 \cite{Poincare92}
($O=\C,\ S_1=[-1,1],\ S_2=(-\infty,-1]\cup [1,+\infty)$ which is, of course,
equivalent to $S_1=(-\infty,0]$, $S_2=[0,\infty)$). The great author was
motivated by his dispute with Borel concerning possible generalizations of the
classical notion of analytic continuation. It is not quite clear who was right
(see the discussion in \cite[Chapter IV, section  21-22]{Valiron54}),  but the
explicit (non-linear) construction of $f_1,f_2$ given in \cite{Poincare92}
(exposed also in \cite{Valiron54}) is very elegant (as was shown in
\cite{MitHen71}, in the Poincar\'e situation there is no linear operator $f\to
(f_1,f_2)$ (from $Hol(O\setminus S)\to Hol(O\setminus S_1)\times Hol(O\setminus
S_2)$ where $Hol(G)$ stands for the space of all functions analytic in an open
set $G$ with the usual topology). Various aspects of the separation of
singularities in the spirit of the Poincar\'e-Aronszajn theorem
(Theorem~\ref{Aron}) are treated in  \cite{Valiron54}, \cite{BerGay91},
\cite{Havin58}, \cite{Havin84} and \cite{Gauthier98}. 

The present article deals with a quantitative aspect of Theorem~\ref{Aron}
related to spaces $H^\infty(G)$ of functions \emph{bounded} and analytic in an
open set $G\subset \C$. Namely we are interested in the case when $f$ in
Theorem~\ref{Aron} is bounded (i.e. $f\in H^\infty(O\setminus S)$) and ask
whether $f_1,f_2$ can be made bounded as well (i.e. $f_j\in H^\infty(O\setminus
S_j),\ j=1,2$).
\begin{definition}
Let $O, S_1,S_2,S$ be as in Theorem~\ref{Aron}. We say $(S_1,S_2)$ is a bounded
separation pair (bs-pair) in $O$ if any $f\in H^\infty(O\setminus S)$ is
representable by the formula
\begin{equation}\label{decomposition}
f=f_1+f_2\text{ in }O\setminus S
\end{equation}
where $f_j\in H^\infty(O\setminus S_j)$.
\end{definition}

The problem which is implicit in the definition can be restated as follows
(just putting $G_j=O\setminus S_j$, $G=O\setminus S$): \emph{given open sets
$G_1,G_2\subset \C$ is it possible to represent an arbitrary $f\in H^\infty(G)$
by \eqref{decomposition} with $f_j\in H^{\infty}(G_j),\ j=1,2$?} 

We want more or less efficient (preferably geometric) conditions imposed on
$(S_1,S_2)$, resp $(G_1,G_2)$ and ensuring \eqref{decomposition} for any $f\in
H^\infty(O\setminus S)$ with $f_j\in H^\infty(O\setminus S_j)$ (resp. $f\in
H^\infty(G_1\cap G_2),\ f_j\in H^\infty(G_j)$).

One of the first results of this kind is due to Polyakov \cite{Polyakov83};  it
was ancillary in his work on bounded extensions of functions $f\in H^\infty(C)$
to the polydisk $\D^n$ (where $\D=\{|z|<1\}$) for an analytic curve $C\subset
\D^n$; see also \cite{HenPol84}. For $G_1=\D\cap \{\Im z >0\}$, $G_2=\D\cap
\{\Re z>0\}$ and any $f\in H^\infty(G_1\cap G_2)$ Polyakov proved
\eqref{decomposition} in $G_1\cap G_2$ where $f_j$ is analytic in $G_j$ and
bounded in $G_j\cap \frac 12 \D$.

In \cite{HavNer01} (whose authors were unaware of \cite{Polyakov83}) and in
\cite{Havin04}  the
bs-pairs were investigated in a systematic way and in much more general
setting. As in \cite{Polyakov83}, in \cite{HavNer01} it was essential that
$S_1,S_2$ meet transversally (a generic situation is presented below on
Figure~\ref{fig3.2}  (section~\ref{sect3.2} below); the results of
\cite{HavNer01} apply to more general transversal configurations). 
In \cite{Havin04} a class of tangent bs-pairs was described. Turning to
the present article we start with some remarks on bs-pairs.

If $S_1,S_2\subset\C$ are compact and disjoint, then $(S_1,S_2)$ is a bs-pair
in $\C$: this is an easy consequence of the Cauchy integral formula. Simple
examples of pairs which are not bs in $\C$ (in particular the Poincare pair
$S_1=(-\infty,0]$, $S_2=[0,\infty)$) are in \cite{HavNer01}.

The bs-problem is related to the Alice Roth Fusion Lemma (see \cite{Gaier87},
\cite{Gamelin69}, \cite{Duady66}, \cite{HavNer01}). In particular, the usual
localization technique (the Vitushkin operator) gives the following result (see
\cite{HavNer01}): \emph{$(S_1,S_2)$ is  a bs-pair in a bounded open set $O$ if
$\dist(S_1,S_2)>0$.}

Dealing with the bs-problem for a pair $(S_1,S_2)$ we may always assume that
$S_1\cap S_2=\emptyset$. (Indeed, $(S_1,S_2)$ is a bs-pair in $O$ iff
$(S_1\setminus S,S_2\setminus s)$ is a bs-pair in $O\setminus s$ where
$s=S_1\cap S_2$.)

A really subtle and interesting situation arises when $S_1\cap S_2=\emptyset$
and $\dist(S_1,S_2)=0$. Let us denote the closure of a set $E$ by $\Clos E$,
and put $K_j=\Clos S_j$, $j=1,2$, $k=K_1\cap K_2$. It is shown in
\cite{HavNer01} that  $(S_1,S_2)$ \emph{is a bs-pair in $O$ iff $(S_1\cap v,
S_2\cap v)$ is a bs-pair in $O$ for a neighbourhood $v$ with respect to $\C$ of
$k$}. Thus only ``the germs of $S_1$ and $S_2$ at $k$'' are responsible for the
bs-property of   $(S_1,S_2)$ in $O$.

In the present paper (as in \cite{HavNer01} and \cite{Havin04}) we concentrate
on the simplest
case when $k$ is a singleton.

We continue the study of ``transversal'' pairs $(S_1,S_2)$ meeting at a point
and obtain a very general result applicable to \emph{arbitrary} relatively
closed parts $S_j$ of the upper half-plane $\C^+$ separated by two rays
$\eta=k\zeta$, $\eta=k_1\zeta$, $\zeta>0$, $0<\Arg(k)<\Arg(k_1)<\pi$. In
\cite{HavNer01} the transversally meeting sets $S_1$ and $S_2$ were supposed to
be parts of certain rectifiable arcs. We believe, however, that the concrete
bounded splittings $f\to (f_1,f_2)$ ($f\in H^\infty(O\setminus S),\ f_j\in
H^\infty(O\setminus S_j)$) constructed in \cite{HavNer01} are of independent
interest; they are also applicable to many \emph{tangent} pairs $(S_1,S_2)$,
see \cite{Havin04}.

Unlike \cite{HavNer01}, we do not exclude the tangency of $S_1$ and $S_2$ and
find rather sharp descriptions of tangent bs-pairs $(S_1,S_2)$ in $\C^+$ when
$k=\{0\}$; in some cases these descriptions yield necessary and sufficient
conditions for a tangent pair of arcs to be bs (Theorem~\ref{teor8}).

We mainly put $O=\C^+$; note that $(S_1,S_2)$ is a bs-pair in $\C^+$ iff
$(\phi(S_1),\phi(S_2))$ is a bs-pair in $\phi(\C^+)$ for a conformal
homeomorphism $\phi$ of $\C^+$.

The article consists of three paragraphs. In \S\ref{p1} we obtain sufficient
conditions for a pair $(S_1,S_2)$ to be bs in $\C^+$ (Theorem~\ref{teor4}). 
Our approach is a reduction to the $\dbar$-problem
\begin{equation}\label{star}
\dbar u = f \dbar \chi
\end{equation}
in an open set $O$ where $f\in H^\infty(O\setminus S)$ is a given function to
be split as in \eqref{decomposition}, and $\chi$ is a cutting factor which is
$\equiv 1$ near $S_1$ and $\equiv 0$ near $S_2$ in $O$; the pair $(S_1,S_2)$ is
bs in $O$ if \eqref{star} has always (for any $f$) a bounded solution $u$.

Theorems~\ref{teor2} and \ref{teor3}  are preparatory; they are corollaries of
deep results by Berndtsson, \cite{Berndtsson92} on bounded solutions of the
$\dbar$-problem in a disc. From them we deduce our main result of \S\ref{p1} 
(Theorem~\ref{teor4}), treating transversal and tangent pairs $S_1,S_2\subset
\C^+$ with $k=\{0\}$; concrete examples are discussed in \ref{sect1.4}. In
\ref{sect1.5bis}, we give a direct and explicit solution of \eqref{star}, not
using Berndtsson's results and thus making our \S\ref{p1} essentially
self-contained. In \ref{sect1.6} we give yet another very explicit proof of
Theorem~\ref{teor4} for pairs $(S_1,S_2)$ in $\C^+$ where $S_1,S_2$ are
separated by an angle with vertex at the origin. In \ref{sect1.7} we briefly
describe quite explicit bounded splittings $f\to(f_1,f_2)$ applicable to a
class of pairs $(S_1,S_2), S_j\subset\C^+$; the splittings of \ref{sect1.5bis},
\ref{sect1.7} are linear in $f\in H^\infty(\C^+\setminus S)$.

The object of \S\ref{p2} are pairs $(S_1,S_2)$ of smooth arcs in an angle $A$
meeting tangentially at its vertex, which are \emph{not} bs in $A$ (and, in
fact, in any domain $O\supset A$).

The construction is a rather involved ``condensation of singularities'', an
accumulation of ``badly splittable'' pairs of arcs in $A$ based  on the Banach
theorem on surjective operators. The results are different depending on whether
the common tangent of $S_1$ and $S_2$ at the vertex is a side of $A$
(Theorem~\ref{teor7}) or is strictly inside $A$ (Theorem~\ref{teor6}). These
theorems show the sharpness of the results of \S\ref{p1}. 

The results of \S\S\ref{p1}--\ref{p2} are combined in \S\ref{p3} to state
Theorem~\ref{teor8} which includes, among other things, the tangential case
when $S_j$ is the graph of a real $\mathcal{C}^{1+\varepsilon}$-function
$\phi_j$ on $[0,b]$ such that
\[
\phi_j(0)=\phi'_j(0)=0,\ 0<\phi_1(t)<\phi_2(t)\quad\text{for }t\in(0,b].
\]
It turns out that $(S_1,S_2)$ is a bs-pair in $\C^+$ iff
\[
\liminf_{x\to 0} \frac{\phi_2(x)-\phi_1(x)}{\phi_1(x)}>0
\]
(i.e. $S_1$ and $S_2$ are hyperbollicaly separated in $\C^+$)

In section~\ref{sect3.2} we show that for a Jordan domain $G=G_1\cap G_2$ where
the domains $G_j$ are Jordan as well, a function $f\in H^\infty(G)$ may exist
which is not representable as $f_1+f_2$ with $f_j\in H^\infty(G_j)$. And in
section~\ref{sect3.3}, we construct an example illustrating the difference
between the case of \emph{connected} sets $S_1,S_2$ (arcs) we dealt with in
Theorem~\ref{teor8}, and \emph{disconnected} $S_j$'s when, apparently, some new
results are needed to grasp the bs-property.

\section{Some classes of bs-pairs in the upper half-plane}\label{p1}
\subsection{Reduction to a $\dbar$-problem}\label{sect1.1}
Let $S_1,S_2$ be disjoint relatively closed subsets of the domain $O$, and
$f\in H^\infty(O\setminus S)$, $S=S_1\cup S_2$. Let $U_1,U_2$ be disjoint open
neighbourhoods of (resp.) $S_1,S_2$ in $O$. Consider a locally Lipschitz
function $\chi$ in $O$ such that
\begin{equation}\label{dollar}
\chi|U_1 =0,\ \chi|U_2=1.
\end{equation}
Extend $f$ to $O$ putting $f|S=0$. Any (distributional) solution $u$ of the
$\dbar$-problem
\begin{equation}
\dbar u = f\dbar \chi\quad \text{in } O
\end{equation}
is a continuous function, since $f\dbar X$ is bounded in $O$. We assume $u\in
\mathcal{C}(O)$ and put 
\[
f_1=f\chi -u,\quad f_2=f(1-\chi) + u\quad \text{in }O,
\]
so that $f=f_1+f_2$ in $O$. Clearly, $f_j$ is continuous in $O\setminus S_j$
and $\dbar f_j=0$ in $O\setminus S_j$, whence it is analytic in $O\setminus
S_j$. It is bounded if $u$ is. We arrive at the following conclusion: 

\emph{if \eqref{star} admits a solution $u\in L^\infty(O)$ for any $f\in
H^\infty(O\setminus S)$, then $(S_1,S_2)$ is a bs-pair with respect to $O$}

\subsection{bs-pairs and a result by Berndtsson}\label{sect1.2}
Bounded solvability of a general $\bar\partial$-equation
\begin{equation}\label{dbar3}
\dbar u =\rho
\end{equation}
in $\D$ was studied in \cite{Berndtsson92}. The result of 
\cite{Berndtsson92} we need can be easily
adapted to $\C^+=\{\Im z>0\}$. By $L^1(E),\ E\subset \C$, we mean $L^1(E,dA)$
where $dA$ is the area element.
\begin{theorem}\label{teor2}
Suppose that $\rho\in L^1(\C^+)$ vanishes in $\C^+\cap \{|z|>R\}$ for a
positive $R$. If 
\begin{itemize}\label{pound}
\item[(a)] $\rho\, dA$ is a Carleson measure in $\C^+$, and
\item[(b)] $\rho \Im z \in L^\infty(\C^+)$,
\end{itemize}
then the (distributional) problem \eqref{dbar3} has a solution $u\in
L^\infty(\C^+)$. 
\end{theorem}
Combining this result with section~\ref{sect1.1} we get
\begin{theorem}\label{teor3}
Suppose $S_1,S_2\subset\C^+$ are bounded, disjoint, and relatively closed.
Suppose there exists a locally Lipschitz real function $\chi$ in $\C^+$
satisfying \eqref{dollar}, vanishing in $\C^+\cap\{|z|>R\}$ and such that
\begin{itemize}
\item[(a)] $|\nabla \chi|\, dA$ is a Carleson measure in $\C^+$,
\item[(b)] $|\nabla \chi(\zeta)|\le C/\Im\zeta$, $(\zeta\in \C^+)$.
\end{itemize}
Then $(S_1,S_2)$ is a bs-pair in any domain $O\subset\C^+$ containing $S$.
\end{theorem}
\begin{proof}
Note that $2|\bar\partial \chi|=|\nabla \chi|$. By Theorem~\ref{teor2} the
$\dbar$-problem~\ref{star}  has a bounded solution in $\C^+$ (thus in $O$), and
section~\ref{sect1.1} applies.
\end{proof}

\subsection{Sets $S_1,S_2$ separated by graphs}\label{sect1.3}
To get a palpable description of some bs-pairs let us assume that there are a
number $\mu>0$ and a Lipschitz nonnegative function $g$ on $\R$ such that
\begin{equation}\label{eq10}
\begin{split}
S_1\subset U_1&=\{\zeta+i\eta \in \C^+:\ \eta<g(\zeta)\},\\
S_2\subset U_2&=\{\zeta+i\eta \in \C^+:\ \eta>(1+\mu)g(\zeta)\}.
\end{split}
\end{equation}
\begin{theorem}\label{teor4}
Suppose $S_1,S_2$ are relatively closed in $\C^+$ and bounded. If \eqref{eq10}
holds then $(S_1,S_2)$ satisfies the assumptions of Theorem~\ref{teor3} and
thus is a bs-pair in any $O\subset\C^+$ containing $S$. 
\end{theorem}
Our assumption \eqref{eq10} means that $S_1,S_2$ are separated by the corridor
$\{\zeta+i\eta\in \C^+:\ g(\zeta)\le\eta\le (1+\mu)g(\zeta) \}$ whose
hyperbolic width is positive (since for $g(\zeta)>0$ the hyperbolic distance
between $\zeta+ig(\zeta)$ and $\zeta+i(1+\mu)g(\zeta)$ exceeds a positive
number not depending on $\zeta\in\R$). We shall see in \S\ref{p2} that
this condition is sharp.
\begin{proof} Put 
\[
\chi_0(\xi,\eta) =
\begin{cases}
0&\text{ if } 0<\eta<g(\xi),\\
\frac{\eta-g(\xi)}{\mu g(\xi)}& \text{ if }\xi\in\R, 0< g(\xi)\le \eta\le
(1+\mu)g(\xi) ,\\
1&\text{ if } \eta>(1+\mu)g(\xi).
\end{cases}
\]
Thus $\chi_0(\xi,\eta)=1$ if $\xi\in\R$, $g(\xi)=0$, $\eta>0$. Clearly
$\chi_0$ is locally Lipschitz on the set $\{\xi+i\eta\in\C^+:\ g(\xi)>0\}$. If
$g(\xi)=0$, $\eta>0$, then $(1+\mu)g(\xi')<\eta'$ for $\xi'\simeq\xi$,
$\eta'\simeq\eta$ whence $\chi_0(\xi',\eta')=1$. We see that $\chi_0\equiv 1$
in a neighbourhood of any point $\xi+i\eta\in\C^+$ with $g(\xi)=0$, and
$\chi_0$ is locally Lipschitz in $\C^+$. By definition $\chi_0|U_1=0$,
$\chi_0|U_2=1$. Take  a large $R>0$ (so that $S\subset\{|\zeta|<R\}$) and then a
$\mathcal C^\infty_0$-function $\chi_1$ which is real and $\equiv 1$ on
$\{|\zeta|<R\}$. Then 
\begin{equation}\label{14}
\chi =\chi_0\chi_1
\end{equation}
satisfies \eqref{dollar}. We only have to check (a) and (b) of
Theorem~\ref{teor3}. But $\nabla
\chi=0$ on $U_1\cup U_2$ and
\[
|\nabla \chi(\xi,\eta)|\le \frac {C(\chi_1)}{\mu g(\xi)} \Bigl( 
\frac{\eta |g'(\xi)|}{g(\xi)}+1 \Bigr)\le C(g,\mu,\chi_1)/\eta
\]
for almost all $\xi\in \{g>0\}$ and $\eta\in (g(\xi),(1+\mu)g(\xi))$ (recall
that $g'\in L^\infty(\R)$). Thus (b), holds. To prove 
(a) consider a Carleson box $B=(a,b)\times (0,b-a)$ ($a,b\in\R,\
a<b$). We have
\[
\int_B |\nabla \chi | dA = \int_{B\cap \{g(\xi)<\eta<(1+\mu)g(\xi)\}}|
\nabla \chi | dA \le C \int_a^b \Bigl(\int_{g(\xi)}^{(1+\mu)g(\xi)}
\frac{d\eta}{\eta}d\xi\Bigr)= C \log(1+\mu)(b-a).
\]
\end{proof}
\subsection{Examples}\label{sect1.4}
\subsubsection{Example 1} For $k>0$ put $g(\xi)=k\xi$ ($\xi\ge 0$), $g(\xi)=0$
($\xi<0$). In this case Theorem~\ref{teor4} is illustrated by pairs $(S_1,S_2)$
separated by \emph{an angle} $\{\xi+ i \eta:\ \xi>0, k\xi <\eta<(1+\mu)k\xi\}$. 
\begin{figure}
\begin{center}
\psfrag{kikprima}{$k'=(1+\mu)k$}
\psfrag{S1}{$S_1$}
\psfrag{S2}{$S_2$}
\psfrag{eta1}{$\eta=k\xi$}
\psfrag{eta2}{$\eta=k'\xi$}
\includegraphics{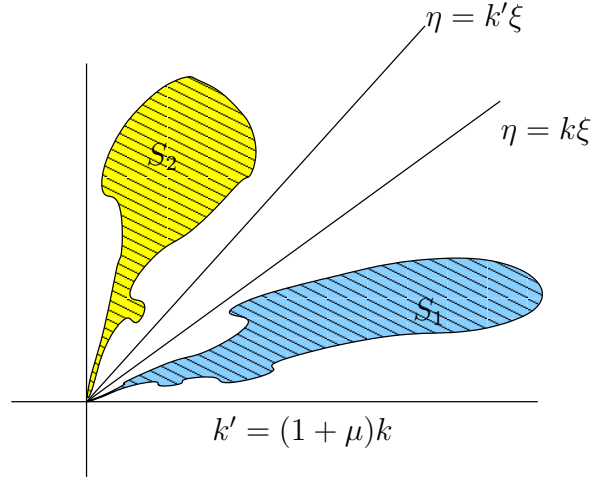}
\end{center}
\caption{Example 1}\label{fig1}
\end{figure}
Compared with \cite{HavNer01} this result is a progress. In \cite{HavNer01} some
transversally meeting bs-pairs were described with an explicit and elementary
splitting formula $f=f_1+f_2$ ($f\in H^{ \infty}(O\setminus S),\ f_j\in
H^\infty(O\setminus S_j)$). But in \cite{HavNer01} certain regularity conditions were
imposed on $S_1,S_2$ (in particular they had to be contained in a union of
rectifiable arcs) whereas Theorem~\ref{teor4} 
allows \emph{arbitrary} relatively
closed sets separated by an angle.
\subsubsection{Example 2} 
\begin{figure}
\begin{center}
\psfrag{S1}{$S_1$}
\psfrag{S2}{$S_2$}
\includegraphics{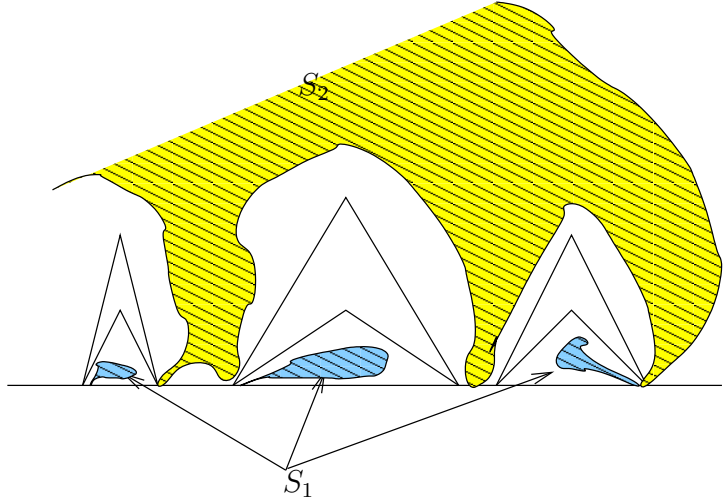}
\end{center}
\caption{Example2}\label{fig2}
\end{figure}

Let $G\subset \R$ bet a bounded open
set. Theorem~\ref{teor4} is applicable to $g(\xi)=\dist(\xi, \R\setminus G)$ 
(see figure~\ref{fig2})
\subsubsection{Example 3} 
\begin{figure}
\begin{center}
\psfrag{eta=g}{$\eta=g(\xi)$}
\psfrag{eta=2g}{$\eta=2g(\xi)$}
\psfrag{S1}{$S_1$}
\psfrag{S2}{$S_2$}
\includegraphics{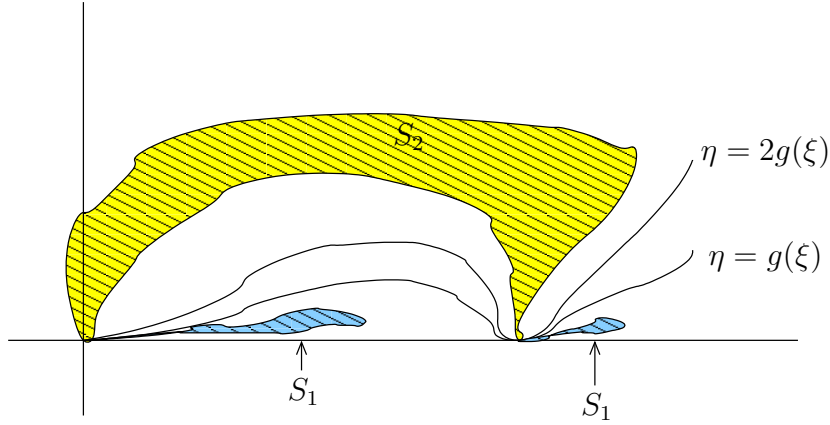}
\caption{Example3}\label{fig3}
\end{center}
\end{figure}

Theorem~\ref{teor4} is also applicable to some
\emph{tangentially} meeting bs-pairs $(S_1,S_2)$. Suppose $g\in
\mathcal{C}^1(\R), g\ge 0, g(0)=g'(0)=0$. On figure~\ref{fig3} $S_1$ and $S_2$
meet tangentially at the origin, but form a bs-pair in $\C^+$ (and in any
subdomain of $\C^+$ containing $S=S_1\cup S_2$). 

\subsection{Theorem~\ref{teor4} revisited: an explicit solution}
\label{sect1.5bis} 
The $\dbar$-estimates that we needed included
uniform bounds in the whole set $\C^+$ and not only on the boundary. That is the
reason that we could not use Carleson estimates, but had to appeal to the more
sophisticated Theorem~\ref{teor4} and we had to assume some extra regularity on
the data $\rho$ in \eqref{dbar3} apart from the Carleson condition. The
estimates
by Berndtsson rely on some a priori estimates and therefore the solution is
not explicit. Nevertheless Jones has found an explicit non-linear 
formula \cite{Jones83} to solve the $\dbar$ equation that gives bounded boundary
values for the solution when the data is a Carleson measure. Seip has
observed (see \cite[Chapter1]{Seip04}) that this formula can be adapted to get
Theorem~\ref{teor4}. For the sake of completeness we present the adapted
Jones solution. 
\begin{proof}
We assume that $|\rho| dA$ is a Carleson measure in 
$\C^+$ and that $\rho \Im z \in L^\infty(\C^+)$. Set
\[
u(z)= \frac{2i}{\pi} \int_{\Im \zeta>0} \frac{1}{z-\zeta}
\frac{\Im \zeta}{z-\bar\zeta}
\exp\Bigl\{\alpha\int_{\Im w\le \Im \zeta}\bigl(\frac{-i}{z-\bar w}+
\frac i{\zeta-\bar w} \bigr)|\rho(w)|dA(w)\Bigr\}\rho(\zeta)dA(\zeta).
\]
The real number $\alpha$ can be chosen freely. 
It is clear that $\dbar u = \rho$ because $u$ if of the form
\[
u(z)=\frac 1\pi \int
\frac{h(z,\zeta)}{z-\bar\zeta}\rho(\zeta)dA(\zeta),
\] 
and for any $\zeta$, $h(\cdot,\zeta)$ is holomorphic and in the
diagonal $h(z,z)=1$.
The hard part is to obtain estimates for $u$.
We take
\[
\alpha^{-1}= \sup_{\Im \zeta >0}\int_{\Im w \le \Im \zeta}
\frac {2\Im \zeta}{|w-\bar \zeta|^2} 
|\rho(w)|dA(w).
\]
This supremum is finite because we assume that $|\rho|dA$ is 
a Carleson measure.
We get
\[
\begin{split}
&|u(z)| \le \\
&\frac{2}{\pi\alpha} \int_{\Im \zeta>0}\! \frac{\alpha}{|z-\zeta|}
\frac{\Im \zeta}{|z-\bar \zeta|} \exp{\Bigl\{\alpha\int_{\Im w \le \Im \zeta} 
\frac{-(\Im w+\Im z)}{|z-\bar w|^2} +\frac{2\Im \zeta}{|\zeta-\bar w|^2}|\rho(w)|dA(w)
\Bigr\}|\rho(\zeta)|dA(\zeta)}.
\end{split}
\]
We use the definition of $\alpha$ and we get
\[
|u(z)| \le \frac{2e}{\pi\alpha} \int_{\Im \zeta>0}\! \frac{\alpha}{|z-\zeta|}
\frac{\Im \zeta}{|z-\bar \zeta|} \exp{\Bigl\{\alpha\int_{\Im w \le \Im \zeta} 
\frac{-(\Im(w)+\Im(z))}{|z-\bar w|^2}|\rho(w)|dA(w)
\Bigr\}|\rho(\zeta)|dA(\zeta)}.
\]
We then write
\[
|u(z)|\le I_1 + I_2,
\]
where
\[
I_1= \frac{2e}{\alpha \pi} \int_{2|z-\zeta| \ge \Im z}\frac{\alpha}{|z-\zeta|}
\frac{\Im \zeta}{|z-\bar \zeta|} \exp{\Bigl\{\alpha\int_{\Im w \le \Im \zeta} 
\frac{-(\Im(w)+\Im(z))}{|z-\bar w|^2}|\rho(w)|dA(w)
\Bigr\}|\rho(\zeta)|dA(\zeta)}
\]
and 
\[
I_2= \frac{2e}{\alpha \pi} \int_{2|z-\zeta| < \Im z}\!\frac{\alpha}{|z-\zeta|}
\frac{\Im \zeta}{|z-\bar \zeta|} \exp{\Bigl\{\alpha\int_{\Im w \le \Im \zeta} 
\frac{-(\Im(w)+\Im(z))}{|z-\bar w|^2}|\rho(w)|dA(w)
\Bigr\}|\rho(\zeta)|dA(\zeta)}.
\]
When $2|z-\zeta|\ge \Im z$, then $|z-\zeta|\ge |z-\bar\zeta|/5$ and thus $I_1$
is bounded by a constant times
\begin{equation}\label{jonestrick}
\frac 1{\alpha} \int_{\Im \zeta>0}\frac {\alpha \Im \zeta}{|z-\bar\zeta|^2}
\exp{\Bigl\{\alpha\int_{\Im w \le \Im \zeta} 
\frac{-\Im(w)}{|z-\bar w|^2}|\rho(w)|dA(w)
\Bigr\}|\rho(\zeta)|dA(\zeta)}.
\end{equation}
If $f$ is any positive function defined on $t\ge 0$ and $R_0=\int_0^
\infty f(t)dt$ (which may eventually be $\infty$), then 
\[
\int_0^\infty f(r) \exp\Bigl(-\int_0^r f(t)\,dt\Bigr)\,
dr=\int_0^{R_0} e^{-t}\, dt\le 1.
\]
We apply this inequality to the function
\[
f(t)=\int_{\R} \frac{\alpha t}{|z-(x-it)|^2}|\rho(x+it)|\, dx,
\]
and we get that \eqref{jonestrick} is bounded by $\alpha^{-1}$.

On the other hand $I_2$ is easy to estimate if we drop the exponential
term and use the fact that 
$|\rho(z)|\le C /\Im z$. 
\end{proof}

\subsection{An explicit solution of \eqref{star} for transversally meeting sets
$S_1,S_2\subset\C^+$}\label{sect1.6}
Here we return to Example~1 of section~\ref{sect1.4} and give another very
simple solution of the problem
\begin{equation}\label{800}
\dbar u=\rho,\ \rho=f\bar\partial \chi,
\end{equation}
where $\chi$ (see \eqref{14}) corresponds to $g$ which is zero on $(-\infty,0]$
and is linear on $[0,\infty)$: $g(\xi)=k\xi$, $\xi\ge 0$, 
for a $k>0$. The sets
$S_1,S_2$ are separated by the angle
\[
A_{k,\mu}=\{(\xi,\eta): k\xi<\eta<(1+\mu)k\xi\},
\] 
and $\partial\chi$ is supported by a sector $T=A_{k,\mu}\cap \{|z|<R\}$. Put
$u=C^\rho -a$, where $C^\rho=\frac 1\pi \rho\star \frac 1z$ is the standard
solution
of \eqref{800} and
\[
a(\zeta)=\frac 1\pi \int_{\C^+} \frac{\bar z}z \frac{\rho(z)}{\zeta-\bar z}
dA(z),\qquad \zeta\in\C^+.
\]
Note that $\rho\equiv0$ off $T$, and $|\rho(z)|\le C/\Im z\simeq 1/|z|$ for
$z\in T$ whence $\rho\in L^1(\C^+,dA)$, and functions $C^\rho, a$ make sense.

Clearly, $\dbar u =\rho$, since $a$ is analytic in $\C^+$, and we only have to
prove that $u$ is bounded in $\C^+$. We have 
\[
u(\zeta)=\frac{2i\zeta}\pi \int_T \frac{\Im z\dbar \chi (z) f(z)
dA(z)}{z(\zeta-z)(\zeta-\bar z)},
\]
whence 
\begin{equation}\label{801}
|u(\zeta)|\le c|\zeta| \int_T \frac{dA(z)}{|z||\zeta-z||\zeta-\bar
z|}=C|\zeta|J(\zeta),\quad \zeta\in \C^+.
\end{equation}
Fix a small $q>0$ (to be specified later). Then 
\[
J(\zeta)\le \int_{T\cap \{|z-\zeta|<q|\zeta|\}}+
\int_{T\cap \{|z-\zeta|>q|\zeta|,\ |z|>q|\zeta|\}}+
\int_{T\cap \{|z-\zeta|>q|\zeta|,\ |z|<q|\zeta|\}}=I+II+III,
\]
Now, $|z|\simeq \Im z$ for $z\in T$. Estimating I we may write $|z|\ge
|\zeta|-|z-\zeta|\ge (1-q)|\zeta|$, 
\[
|\bar z-\zeta|\ge 2 \Im z -|z-\zeta|\ge c|z|-q|\zeta|\ge c'|\zeta|
\]
where $c=c(k,\mu)>0$, $c'=c(1-q)-q>0$ if $q$ is small, and 
\begin{equation}\label{802}
I\le \frac{c_1}{|\zeta|^2}\int_{ |z-\zeta|<q|\zeta|}
\frac{dA(z)}{|\zeta-z|}=\frac{c_2}{|\zeta|}.
\end{equation}
The H\"older inequality (with $p=3$, $q=3/2$) and $|\bar z-\zeta|\ge|z-\zeta|$
give 
\begin{equation}\label{803}
II\le \int_{ |z-\zeta|>q|\zeta|,\ |z|>q|\zeta|}\frac {dA(z)}
{|z||z-\zeta|^2}\le 
\left(\int_{|z|>q|\zeta|} \frac {dA(z)}{|z|^3} \right)^{1/3}
\left(\int_{|z-\zeta|>q|\zeta|}\frac
{dA(z)}{|z-\zeta|^3}\right)^{2/3}=\frac{c_3}{|\zeta|}.
\end{equation}
At last, 
\begin{equation}\label{804}
III\le \int_{|z-\zeta|>q|\zeta|,\ |z|<q|\zeta|}\frac {dA(z)}
{|z||z-\zeta|^2}\le\frac 1{q^2|\zeta|^2}\int_{|z|<q(|\zeta|}\frac{dA(z)}{|z|}
=\frac{c_4}{|\zeta|}. 
\end{equation}
Combining \eqref{801} with \eqref{802}--\eqref{804} we see that $u$ is bounded.
\qed
\subsection{Another explicit and bounded solution to
equation \eqref{star}}\label{sect1.7}
Here we again exploit the special form of the right side of
\eqref{star} (recall that Theorem~\ref{teor3} is aimed at
\emph{general} right sides $\rho\in L^1(\C^+)$, but we only deal with
$\rho=f\dbar \chi$, where $f\in H^\infty(\C^+\setminus S)$ and $\chi$
is a smooth cutting factor).

This time we consider a function $g:\R\to[0,+\infty)$ and assume
$g\in\mathcal C^{1+\varepsilon}(\R)$, $\varepsilon>0$, $g\equiv 0$ on
$(-\infty,0]$, $g(x)>0$ for $x>0$ (in Theorem~\ref{teor4} $g$ was a
Lipschitz function and it could vanish at some positive points). The
functions $\chi_0,\ \chi_1,\ \chi$ are the same as in
Theorem~\ref{teor4}. Fix a small $b>0$ so that $\chi_1\equiv 1$ on the
corridor
\[
G_b=\{\xi+i\eta:\ 0<\xi<b,\ g(\xi)<\eta<(1+\mu)g(\xi)\},
\]
and $\chi_0\equiv \chi$ on $G_b$. For $t\in[1,1+\mu]$ put
$\gamma_t(\xi)=\xi+itg(\xi)$, $0\le\xi\le b$. 

Suppose $f$ is as in Theorem~\ref{teor4} ($f\in L^\infty(\C^+)$,
$f|S=0$, $f$ is analytic in $\C^+\setminus S$). Put 
\[
u(\zeta)=\frac 1\pi \int_{\C^+}\frac{f(z)\dbar\chi(z)}{\zeta-z}dA(z)+
\frac{\mu}{2\pi i}\int_{\bar\gamma_1}\frac{f(\bar z)dz}{\zeta-z},\ \zeta\in\C^+.
\]
\emph{Then $u$ is bounded in $\C^+$ and satisfies \eqref{star}}. 

Recall
that $\gamma_1$ is the path $t\to t+ig(t)$, $t\in[0,b]$, so that the
path $\bar\gamma_1$ (i.e. the reflected graph of $g$) does not
intersect $\C^+$, and the contour integral is analytic in $\C^+$. At
the same time the double integral represents ``the standard solution
of \eqref{star}'' (note that $\dbar\chi$ is summable in $G_b$, see the
last estimate in Section~\ref{sect1.3}). Thus our $u$ staisfies
$\eqref{star}$ in $\C^+$, and the only problem is its boundedness. The
standard solution may blow up at the origin, but its growth is
counterbalanced by the contour integral. A detailed proof will be
given in \cite{HavOrt05}. 

Our definition of $u$ combined with Section~\ref{sect1.1} generates a
\emph{linear} operator $f\to(f_1,f_2)$ splitting $f\in H^
\infty(\C^+\setminus S)$ into the sum of functions $f_j\in
H^\infty(\C^+\setminus S_j)$, $j=1,2$.

\section{Negative results}\label{p2}
In this paragraph we describe some pairs $(S_1,S_2)$ of subsets of an angle $A$
which are \emph{not} bs-pairs in $A$ (Theorems~\ref{teor6}, \ref{teor7} in
section \ref{sect2.4.5} and \ref{sect2.5}). Combined with the results of
\S\ref{p1} they give  some necessary and sufficient geometric conditions for two
smooth
\emph{graphs} to form a bs-pair in $\C^+$ (see \S\ref{p3}). We start with preliminary
technical results on ``badly splittable'' pairs of arcs $K_1,K_2$ in a domain
$O$ which means the existence of $\phi\in H^\infty(O\setminus K)$, $K=K_1\cup
K_2$ with $|\phi|\le 1$ in $O\setminus K$ and such that a representation
$\phi=\phi_1+\phi_2$, $\phi_j\in H^\infty(O\setminus K)$, is only possible with
a very big $\sup \{|\phi_1(\zeta)|:\ \zeta\in O\setminus K_1\}$. If
$S_1,S_2\subset O$ contain the elements $K_1,K_2$ of arbitrarily badly
splittable pairs $(K_1,K_2)$ (i.e. $S_j\supset K_j$, $j=1,2$), then
$(S_1,S_2)$ is not a bs-pair in $O$. These vague considerations are made precise
in sections~\ref{sect2.1}--\ref{sect2.4}, and then applied to quite concrete
Theorems~\ref{teor6}, \ref{teor7}.
\subsection{Cells and their rotundities}\label{sect2.1}
A \emph{cell} is by definition a pair $(g,A)$ where $g$ is a Jordan domain with
a rectifiable boundary $\partial g$, and $A$, ``the center'' of $g$, is a point
in $g$.  Sometimes we write $g$ instead of $(g,A)$. Put
\[
\rho(g)(=\rho((g,A)))=\frac {2\pi\dist(A,\partial g)}{l(\partial g)}
\] 
where $l(\partial g)$ denotes the length of $\partial g$. Clearly, $\rho(g)\le
1$, a geometrically obvious fact (an instant proof is given by the formula 
$1=(2\pi i)^{-1}\int_{\partial g}\frac{dz}{z-A}$); the equality $\rho(g)=1$
occurs iff $g$ is a disc centered at $A$. We call $\rho(g)$ \emph{the 
rotundity} of the cell $g(=(g,A))$.

\subsection{Functions at a large distance from $H^\infty(g)$ in a cell
$g$}\label{sect2.2}
Let $K$ be a compact subset of $g$ containing $A$. We assume $A$ to be a
boundary point of $K$ (in our applications $K$ will be just a simple arc)

\begin{lemma}\label{lema2.1}
Let $\phi$ be a function analytic in $\C\setminus K$, $\phi(\infty)=0$. Then any
$h\in H^\infty(g)$ satisfies 
\begin{equation}\label{instar}
\|\phi -h\|_{g\setminus K} \ge \frac{\rho(g)}{2} \limsup_A |\phi|,
\end{equation}
where $\|\cdot\|_E$ denotes $\|\cdot\|_{\infty,E}$.
\end{lemma}
In other words, if $\phi$ is very big near $A$ and $(g,A)$ is sufficiently
rotund, then \emph{any} $h\in H^\infty(g)$ is far away from $\phi$ in
$H^\infty(g\setminus K)$. 
\begin{proof}
Any $h\in H^\infty(g)$ has the angular boundary value $h(g)$ at almost every
$q\in\partial g$, and
\[
h(p)\equiv \frac 1{2\pi i} \oint_{\partial g} \frac{h(q)}{q-p}\,
dq\qquad\text{for }p\in g.
\]
Moreover,  $\oint_{\partial g} \frac{\phi(q)}{q-p}\equiv 0$ for $p\in g$.
Thus putting $\|\cdot\|=\|\cdot\|_{g\setminus K}$, we have
\[
\begin{split}
|\phi(p)|\le |h(p)| + \|\phi-h\|=\\
\frac 1{2\pi} \Bigl|\oint_{\partial g} \frac{h(q)-\phi(q)}{q-p}dq\Bigr|+
\|\phi-h\|\le \\
\bigl(\frac{1}{\rho(g)}+1\bigr)\|\phi-h\|\le \frac 2{\rho(g)}\|\phi-h\|.
\end{split}
\]
Letting $p\to A$ we get \eqref{instar}.
\end{proof}
\subsection{Pairs not admitting bounded separation in a domain (an abstract
scheme)}\label{sect2.3}
Lemma~\ref{lema2.1} suggests a method to construct pairs $(S_1,S_2)$ which are
not bs in a domain $O$.
\subsubsection{}\label{3star} Suppose $S_1,S_2$ are relatively closed disjoint and nowhere
dense parts of a domain $O$. Suppose 
there is a number $C$ such that for any (big) $M>0$ there exists a cell
$g\bigl(=(g,A)\bigr)$, $g\subset O$, and a pair $(K_1,K_2)$ of compacts such
that $K_j\subset S_j$, $j=1,2$, and $K_1\subset g$, $\rho(g)\ge C^{-1}$, and for
a pair $(\phi_1,\phi_2)$ of functions analytic, respectively, in $\C\setminus
K_j$, $j=1,2$, we have
\[
\phi_1(\infty)=\phi_2(\infty)=0,\ |\phi_1+\phi_2|\le C \text{ in }
\C\setminus K,\ K=K_1\cup K_2,
\]
whereas
\begin{equation}\label{4star}
\limsup_A |\phi_1|>M.
\end{equation}
\begin{lemma}\label{lema2.2} If $S_1,S_2$ and $O$ enjoy property ~\ref{3star},
then $(S_1,S_2)$ is not a bs-pair in $O$.
\end{lemma}
Note that $\phi_j\in H^\infty(\C\setminus K_j)$, $j=1,2$, since $K_1,K_2$ are
disjoint.
\begin{proof}
Put $S=S_1\cup S_2$. If $(S_1,S_2)$ were a bs-pair in $O$, then the operator
$(f_1,f_2)\to (f_1+f_2)(O\setminus S)$ from $H^\infty(O\setminus S_1)
\times H^\infty(O\setminus S_2))$ to $H^\infty(O\setminus S)$ would be 
surjective, and, by the Banach theorem, we could 
find a number $N>0$ such that any $f\in
H^\infty(O\setminus S)$ would admit a splitting
\begin{equation}\label{!!}
f=f_1+f_2\qquad\text{in }O\setminus S, 
\end{equation} 
with $f_j\in H^\infty(O\setminus S_j)$, $\|f_j\|_{O\setminus S_j}\le N
\|f\|_{O\setminus S}$.
 
But \ref{3star} forbids the existence of such number $N$. Indeed, let $K_j\
\phi_j$ be as in \ref{3star} where a big $M=M(N)$ will be specified later.
Put $f_j^0=\phi_j|(O\setminus K_j)$, $j=1,2$, and $f=f_1^0+f_2^0$ in $O\setminus
K$, $K=K_1\cup K_2$. Then for any representation \eqref{!!} of $f$ we get 
\[
\begin{split}
f_1&=f_1^0-h \text{ in } O\setminus S_1\\
f_2&=f_2^0+h \text{ in } O\setminus S_2
\end{split}
\]
where $h\in H^\infty(O)$, since $f_1^0-f_1$ and $f_2^0-f_2$ are mutual analytic
continuations from $O\setminus S$ across $S_2$ and $S_1$. But $f_1$ is in fact
defined and analytic not only on $O\setminus S_1$, but in $O\setminus K_1$,
since $f_1^0\in H^\infty(O\setminus K_1)$, $h\in H^\infty(O)$. Moreover, 
$\|f_1\|_{O\setminus K_1}= \|f_1\|_{O\setminus S_1}$, since $S_1$ has no
interior points. Now since $h\in H^\infty(g)$ and $\rho(g)\ge C^{-1}$, by
Lemma~\ref{lema2.1} and \eqref{4star} we have 
\[
\|f_1\|_{O\setminus S_1}=\|f_1\|_{O\setminus K_1}\ge \frac {M}{2C}.
\]
For an arbitrary $N>0$ take $M>2NC^2$. We get a contradiction: given $N>0$, no
splitting \eqref{!!} of $f$ with $f_j\in H^\infty(O\setminus S_j)$ satisfies 
$\|f_1\|_{O\setminus S_1}\le N \|f\|_{O\setminus S}$ (recall that
$\|f\|_{O\setminus S}\le C$).
\end{proof}
\subsubsection{}\label{sect2.3.2}
Now we are ready to describe a general scheme of constructing triples
$(O,S_1,S_2)$ where $O$ is a domain, $S_j$ its relatively closed and disjoint
subsets such that $(S_1,S_2)$ is not a bs-pair in $O$. Namely, suppose we have
constructed two families $(K_1^x)_{x\in E}$, $(K_2^x)_{x\in E}$ of compact
simple arcs and two families $(\phi_1^x)_{x\in E}$, $(\phi_2^x)_{x\in E}$ of
functions in, respectively, $H^\infty(\C\setminus K_j^x)$, $j=1,2$, and a family
of points $(A^x)_{x\in E}$ such that
\begin{itemize}
\item[(a)] $K_1^x\cap K_2^x=\emptyset$;
\item[(b)] $\phi_j^x(\infty)=0,\ j=1,2$;
\item[(c)]\label{(c)} $\sup_{x\in E}\|\phi^x_1+\phi_2^x\|_{\C \setminus K^x}<+\infty$,
$K^x:=K_1^x\cup K_2^x$;
\item[(d)] $A^x\in K_1^x$ and $\sup_{x\in E}\limsup_{A^x}|\phi_1^x|=+\infty$.
\end{itemize}
The last element of our construction is a family $(g^x,A^x)_{x\in E}$ of cells
such that
\begin{itemize}
\item[(e)] $K_1^x\subset g^x$, $x\in E$, and $g^x$ are uniformly
rotund (i.e. $\inf_{x\in E} \rho(g^x)>0$).
\end{itemize}
Any family of quintuples $(K_1^x,K_2^x,\phi_1^x,\phi_2^x,g^x)$ enjoying 
(a)--(e) generates plenty of ``bad'' triples $(O,S_1,S_2)$. Indeed, for
\emph{any} domain $O$ containing $\cup_{x\in E}(g^x\cup K_2^x)$, any pair
$(S_1,S_2)$ of its disjoint relatively closed and nowhere dense subsets is not a
bs-pair with respect to $O$ provided 
\[
S_1\supset \cup_{x\in E}K^x_1,\qquad S_2\supset \cup_{x\in E}K^x_2.
\]
This is an immediate consequence of Lemma~\ref{lema2.2}.

Sections~\ref{sect2.1}--\ref{sect2.3} are an abridged version of \S\ref{p1} in
\cite{Havin04}. But the technique of section~\ref{sect2.4} is quite different
from \cite{Havin04} and results in removing a logarithmic factor from Theorem~4
of \cite{Havin04}; our Theorem~\ref{teor8} is also based on this improvement.

\subsection{A realization of the general scheme: pairs of graphs}
\label{sect2.4}
In what follows the arcs $K_j^x$ will be pieces of the graphs of real functions
$\phi_j,\ j=1,2$ defined on $[0,b]$. We assume
\begin{equation}\label{16}
\varphi_j\in\mathcal{C}^{1+\varepsilon}([0,b]),\ \varphi_j(0)=\varphi_j'(0)=0,\
j=1,2,\ \varphi_1(t)<\varphi_2(t)\ \text{ for }t\in (0,b].
\end{equation}
It will be also convenient to assume
\begin{equation}\label{17}
|\varphi_j'(t)|<1/2,\ t\in[0,b] 
\end{equation}
(this assumption doesn't affect the generality).
\begin{figure}
\begin{center}
\psfrag{K2x}{$K_2^x$}
\psfrag{K1x}{$K_1^x$}
\psfrag{0}{$0$}
\psfrag{x}{$x$}
\psfrag{x1}{$x_1$}
\psfrag{x2}{$x_2$}
\psfrag{X}{$X$}
\psfrag{b}{$b$}
\includegraphics{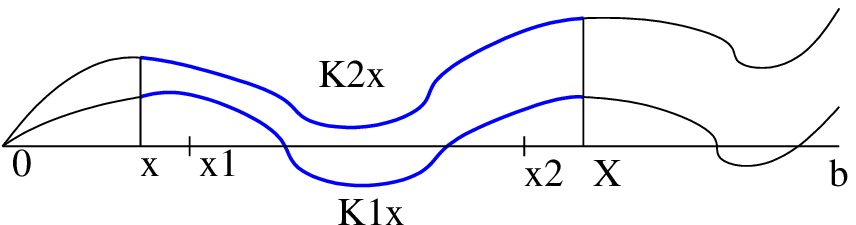}
\end{center}
\caption{$x_1-x=X-x_2=h,\ X=X(x),\ x\to 0$}\label{fig2.1}
\end{figure}

\subsubsection{}\label{sect2.4.1} 
We put
\begin{equation}\label{18}
K_j^x = \{t+i\varphi_j(t): x\le t\le X\},\ x\in (0,b]
\end{equation}
The choice of $X=X(x)$ will be specified later.

\begin{figure}
\begin{center}
\psfrag{0}{$0$}
\psfrag{x}{$x$}
\psfrag{x1}{$x_1$}
\psfrag{x2}{$x_2$}
\psfrag{X}{$X$}
\psfrag{b}{$b$}
\psfrag{1}{$1$}
\includegraphics{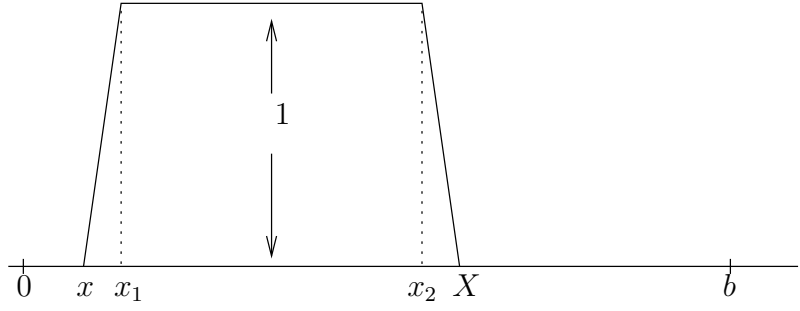}
\end{center}
\caption{$x_1-x=X-x_2=h(x)$}\label{fig2.2}
\end{figure}
To define the functions $\phi_j^x$ (see (c) on page~\pageref{(c)}) we need a
piecewise linear function $f^x$ defined on $\R$ as shown on Figure~\ref{fig2.2},
so that $f^x\equiv 0$ on $\R\setminus(x,X)$, $f^x\equiv 1$ on $[x_1,x_2]$, $f^x$
is linear on $[x,x_1]$, $[x_2,X]$; the positive number $h=h(x)$ will be chosen
later. Now put 
\[
\tilde{f}_j^x(\zeta)=f^x(\Re \zeta)\qquad (\zeta\in \C).
\]
Define the functions $\phi_j^x$ by the formulas
\begin{equation}\label{19}
\begin{split}
\phi_j^x(\zeta)=&(-1)^{j-1}\frac 1{2\pi i} \int_{K_j}
\frac{\tilde{f}_j^x(z)}{z-\zeta}dz=\\
&(-1)^{j-1}\frac 1{2\pi i} \int_x^X \frac{f^x(t)(1+i\varphi_j'(t))}
{(t-\xi)+i(\varphi_j(t)-\eta)}\, dt,\quad \zeta=\xi+i\eta\in \C\setminus K^x_j.
\end{split}
\end{equation}
The first (contour) integral in \eqref{19} is oriented from the left to the
right.

In fact the integral \eqref{19} can be easily evaluated, so that $\phi_j^x$
become explicit linear combinations of some elementary functions of $\zeta$
(compositions of $\zeta\log\zeta$ and $\log\zeta$ with some M\"obius functions).
But we prefer to keep the integral representation of $\phi_j^x$, the estimates
of their explicit versions not being any shorter. Getting rid of integration
makes it possible to apply the scheme of Section~\ref{sect2.3} to general
continua $S_1$, $S_2$ (not just curves), see \cite{Havin04}; but it also makes
the description of bs-pairs of \emph{graphs} not as precise as our
Theorem~\ref{teor8}. 

Put $A^x(=A)=x+i\varphi_1(x)$. Functions $\phi_j^x$ belong to $H^\infty(\C\setminus
K_j^x)$ and vanish at infinity; $\phi_j^x$ is continuous at both sides of the
arc $K_j^x$. In particular, $\lim_{A_x}\phi_1^x$ exists and coincides with 
$\frac 1{2\pi i}\int_{K_1^x}\frac{\tilde{f}_1^x(z)dz}{z-A^x}$. We denote it by
$\phi_1(A)$.

We are going to show that $|\phi_1(A)|$ is big if $h<< X-x$. Indeed, 
\[
\phi_1(A) = \Bigl( \int_x^{x_1}+\int_{x_1}^{x_2}+\int_{x_2}^X \Bigr)
\frac{f(t)}{2\pi i} \frac{(1+i\varphi_1'(t))dt}{(t-x)+i(\varphi_1(t)-\varphi_1(x))}= I +
II+III,
\]
\[
|I| \le \frac 1{2\pi} \int_x^{x+h} \frac{t-x}h \frac{3/2 dt}{t-x}<1.
\]
To estimate III note that any $t\in
(x_2,X)$ satisfies $X-t<t-x$, and
\[
|III| \le \frac 1{2\pi} \int_{X-h}^X \frac{X-t}h \frac{3/2 dt}{t-x} <1.
\]
Now, 
\[
|II|= \frac 1{2\pi} \Bigl| \int_{x_1}^{x_2}
\frac{d(t+i\varphi_1(t))}{(t+i\varphi_1(t))-A}
\Bigr|\ge \frac 1{2\pi}\Bigl|\log 
\frac{|(x_2-x)+i(\varphi_1(x_2)-\varphi_1(x))|}{|(x_1-x)+i(\varphi_1(x_1)-\varphi_1(x))|}
\Bigr|.
\]
The numerator of the last fraction is greater or equal than $x_2-x$, the
denominator is less or equal than $3/2(x_1-x)$ whence
$|II|\ge \frac 1{2\pi} \log \frac{x_2-x_1}{2h}$. But
$h/(x_2-x_1)=(h/(X-x))(1-2h/(X-x))^{-1}$, so that
\begin{equation}\label{20a}
\lim_{x\to 0} |\phi_1^x(A^x)|=+\infty\text{ if } \lim_{x\to 0}
\frac{h(x)}{x-X(x)}=0.
\end{equation}
\subsection{Upper estimates of $|\phi_1^x+\phi_2^x|$}\label{sect2.4.2}
Put 
\[
K^x=K_1^x\cup K_2^x,\ \phi^x=\phi_1^x+\phi_2^x.
\]
(Sometimes we omit the index $x$). We now have to concentrate on the upper
estimate of $|\phi|$ in $\C\setminus K$ (see (c) on page \pageref{(c)}). We may
apply the elementary maximum modulus principle to $\phi$, since $\phi$ is
continuous up to the sides of the arcs $K_1,K_2$. Expressing the boundary
values of $\phi$ on $K$ by the Sohotsky-Plemelj 
formulas 
and using the estimate $|f|\le 1$, we reduce our problem to the upper
estimate of $|J(z_0)|$ where
\[
J(z_0)=p.v. \int_K\tilde f(z)\frac{dz}{z-z_0},\qquad z_0\in K,
\]
$K_1$ is oriented from the left to the right, $K_2$ is oriented in the opposite
direction, so that, hopefully, the contributions of $K_1$ and $K_2$ will be
mutually tempered. Put
\[
z_j(t)= t +i \varphi_j(t),\qquad t\in[0,b],
\]
and suppose 
\[
z_0\in K_1,\ z_0=z_1(t_0),\ t_0\in[0,b].
\]
(the case $z_0 \in K_2$ is symmetric). We have
\[
J(z_0) =p.v.\int_x^X f(t) \K(t,z_0)\, dt,
\]
where
\[
\K(t,z_0)=\frac{z_1'(t)}{z_1(t)-z_0}-\frac{z_2'(t)}{z_2(t)-z_0}.
\]
For a smooth function $\lambda$ of a real variable $t$ put 
\[
R_{t_0} (\lambda)(t)= \lambda(t)-\lambda(t_0)-\lambda'(t)(t-t_0).
\]
Put $\Delta =\varphi_2-\varphi_1$. Using this notation the reader can check the
following identity:
\begin{equation}\label{21}
\K(t,z_0)=\K_1(t,z_0)+\K_2(t,z_0)
\end{equation}
where
\begin{equation}\label{22}
\begin{split}
\K_1(t,z_0)&=\frac{iR(\Delta)+\varphi'_2 R(\varphi_1)-\varphi_1' R(\varphi_2)}
{(z_1(t)-z_0)(z_2(t)-z_0)},\\
&\K_2(t,z_0)=\frac{i\Delta(t_0)z_1'(t)}{(z_1(t)-z_0)(z_2(t)-z_0)},\ R=R_{t_0}.
\end{split}
\end{equation}
An easy estimate of the kernel $K_1$ follows from the inequalities
\[
|z_j(t)-z_0|\ge |t-t_0|,\ |R(\varphi_j)(t)|\le c(\varphi_j)
|t-t_0|^{1+\varepsilon},
\ j=1,2,
\]
since $\varphi_j$ (and $\Delta$) are in $\mathcal{C}^{1+\varepsilon}$; thus,
for
$t_0\in[0,b]$, 
\begin{equation}\label{23}
\Bigl|p.v.\int_x^X f(t) \K_1(t,z_0)\, dt \Bigr|\le c
\int_x^X\frac{dt}{|t-t_0|^{1+\varepsilon}}\le c'
\end{equation}
where $c,c'$ depend only on $\varphi_1$, $\varphi_2$ and $\varepsilon$ (but not
on $z_0$) for any $f$ satisfying $|f(t)|\le 1$ on $[0,b]$.

Turn now to 
\begin{equation}\label{24}
J_2=p.v.\int_x^X \K_2(t,z_0)f(t)\, dt=\int_{[x,X]\setminus a} +\int_a =j_1 +j_2
\end{equation}
where $a=[x,X]\cap [t_0-\Delta(t_0), t_0+\Delta(t_0)]$.
The integral $j_1$ is easy (recall that $|f|\le 1$, $|\varphi_1'|\le 1/2$):
\begin{equation}\label{25}
|j_1|\le 2 \Delta (t_0) \int_{|t-t_0|>\Delta(t_0)} \frac {dt}{(t-t_0)^2}=4.
\end{equation}
But $j_2$, unlike preceding estimates, requires some special properties of $f$
and additional restrictions on the proximity of $K_1$ and $K_2$ (i.e. smallness
restrictions on $\Delta$). To estimate $j_2$ put  $F(t)=f(t)/(z_2(t)-z_0)$; we
get 
\begin{equation}\label{30}
\begin{split}
j_2=&i\Delta(t_0) p.v. \int_a F(t) \frac{z_1'(t)\, dt}{z_1(t)-z_0}=\\
&i\Delta(t_0) \int_a \frac{F(t)-F(t_0)}{z_1(t)-z_1(t_0)}dz_1(t)+i
\Delta(t_0)F(t_0) p.v. \int_{K_1^a} \frac {dz}{z-z_0}= I + II,
\end{split}
\end{equation}
where $K_1^a=\{t+i\varphi_1(t):\ t\in a\}$. Now,
\begin{equation}\label{33}
|I| \le \Delta(t_0) \int_a \left|\frac{F(t)-F(t_0)}{t-t_0}\right|
(1+|\varphi_1'(t)| )\, dt\le 2\Delta(t_0) \|F'\|_{\infty,a} |a|,
\end{equation}
\begin{equation}\label{34}
|F'(t)|\le |f'(t)|/|z_2(t)-z_1(t)| + |f(t)|
\frac{1+|\varphi_2'(t)|}{|z_2(t)-z_1(t_0)|^2}.
\end{equation}
But 
\[
|z_2(t)-z_1(t_0)|\ge |\varphi_2(t)-\varphi_1(t_0)|\ge \Delta(t_0)
-\|\varphi_2'\|_{\infty} |t-t_0|\ge \Delta(t_0)/2,
\]
since $\|\varphi_2'\|<1/2$, $t\in a$. From \eqref{33}, \eqref{34} and the
estimates 
$\|f\|_{\infty}\le 1$, $\|f'\|_{\infty}\le 1/h$, $|a|\le 2\Delta(t_0)$ we
conclude that 
\[
|I| \le 2\Delta(t_0) \left(\frac 2{h\Delta(t_0)}+\frac 8{\Delta(t_0)^2}
\right)2\Delta(t_0)=\frac{8\Delta(t_0)}h
+32\le 40, 
\]
if
\begin{equation}\label{35}
\Delta(t_0)\le h,\quad\text{for any } t_0\in [x,X],
\end{equation}
an important restriction expressing the proximity of $\varphi_1$ and $\varphi_2$
(note
that $h<< X-x$ by \eqref{20a}).

We are left now with the integral II (see \eqref{30}). From
$|F(t_0)|\Delta(t_0)=f(t_0)$ we conclude 
\[
\begin{split}
|II|\le & f(t_0) \left| p.v. \int_{K^a_1}\frac{dz}{z-z_0}
\right|,\\ & a=[x,X]\cap [t_0-\Delta(t_0), t_0+\Delta(t_0)],\ K_1^a=K_1\cap
\{\zeta:\ \Re\zeta\in a\},\ t_0\in [x,X].
\end{split}
\]
\subsubsection{}\label{sect2.4.3}
We have to consider two cases:
\begin{itemize}
\item[Case 1:] $x\le t_0-\Delta(t_0)\le t_0+\Delta(t_0)\le X$;
\item[Case 2:] $t_0-\Delta(t)<x$ or $t_0+\Delta(t_0)>X$.
\end{itemize} 
In Case 1 $a=[t_0-\Delta(t_0), t_0+\Delta(t_0)]$, and 
\[
\begin{split}
|II|\le &\left|p.v.\int_{K_1^a} \frac{dz}{z-z_0}\right|=
\left| \log \frac{z_1(t_0+\Delta(t_0))-z_1(t_0)}
{z_1(t_0-\Delta(t_0))-z_1(t_0)}
\right|\le \\
\le &\left| \log \frac{|z_1(t_0+\Delta(t_0))-z_1(t_0)|}
{|z_1(t_0-\Delta(t_0))-z_1(t_0)|}
\right|+2\pi=\\
=&\frac 12 
\left| \log \frac{(\Delta(t_0))^2+(\varphi_1(t_0+\Delta(t_0))-\varphi_1(t_0))^2}
{(\Delta(t_0))^2+(\varphi_1(t_0-\Delta(t_0))-\varphi_1(t_0))^2}
\right|+2\pi\le B,
\end{split}
\]
$B$ being an absolute constant, since $|\varphi_1(t_0\pm
\Delta(t_0))-\varphi_1(t_0)|\le \Delta(t_0)/2$ whence the fraction under the
logarithm is in $(4/5,5/4)$.

In Case 2 we suppose, e.g., that $t_0-\Delta(t_0)<x$ (the case
$X<t_0+\Delta(t_0)$ is symmetric). Then $0\le t_0-x< \Delta(t_0)\le h$ by
\eqref{35}, so that $t_0\in [x,x_1]$ and $f(t_0)=(t_0-x)/h$. Now
$a=[x,t_0+\Delta(t_0)]$ (note that $t_0+\Delta(t_0)>X$ is impossible, because
this inequality implies $2\Delta(t_0)=(t_0+\Delta(t_0))-(t_0-\Delta(t_0)\ge 
X-x>2h$, a contradiction with \eqref{35}). Thus,
\begin{equation}\label{40}
|II|\le \frac{t_0-x}{h}\left[\left|\log \frac{|z_1(t_0+\Delta(t_0))-z_1(t_0)|}
{|z_1(x)-z_1(t_0)|}
\right|+2\pi \right]
\end{equation}
The numerator of the last fraction is between $\Delta(t_0)$ and $2\Delta(t_0)$
(since $|\varphi_1'|<1/2$), its denominator is between $t_0-x$ and $2(t_0-x)$ so
that, by \eqref{35} and \eqref{40}, 
\[
|II|\le \frac{t_0-x}{\Delta(t_0)}\left[
\log\frac{\Delta(t_0)}{t_0-x}+\log 2 +2\pi
\right]\le e^{-1}+\log 2 +2\pi
\]
(since $y|\log y|\le e^{-1}$ for $y\in (0,1]$).

Summing up, \emph{our functions $\phi_1^x$, $\phi_2^x$ satisfy conditions
(a)--(d) of Section~\ref{sect2.3.2}, if \eqref{20a} and \eqref{35} are
fulfilled}.
\subsubsection{}\label{sect2.4.4}
We now have to specify $h(x)$, $X(x)$ to satisfy \eqref{20a} and \eqref{35}, and
then construct the cells $g^x$ centered at $A^x=x+i\varphi_1(x)$ and satisfying
condition (e) in \ref{sect2.3.2}. These definitions will be given twice, one
time for the angle $A_k=\{\zeta=\xi+i\eta\in\C:\ \xi>0,\ |\eta|<k\xi\}$, and the
second for its upper half $A_k^+=A_k\cap \C^+$.

For any $R>0$ we may assume the graphs of $\varphi_1$ and $\varphi_2$ are in
$A_k$, replacing $b$ by a smaller number if needed. Put
\[
\epsilon(x)=\sup \{\frac{|\varphi_1(t)|+|\varphi_2(t)|}t:\ 0<t\le 2x\},\
X(x)=2x,\ h(x)=2\varepsilon(x)x.
\]
Then $\lim_{x\to 0}\frac{h(x)}{X(x)-x}=\lim_{x\to0} 2\epsilon(x)=0$ (see
\eqref{16} in section~\ref{sect2.4}), 
and we get \eqref{20a}. If $t\in [x,X(x)]$,
then $\Delta(t)=\varphi_2(t)-\varphi_1(t)\le |\varphi_2(t)|+|\varphi_1(t)|
\le \epsilon(x)t\le h(x)$, so that \eqref{35} holds. Put $g^x=A_k\cap \{
\zeta\in\C:\ \Re\zeta<3x\}$. The inclusion $K_1^x\subset g^x$ is obvious if
$x>0$ is small (again by \eqref{16} and \eqref{18} in section~\ref{sect2.4}) as
is the uniform rotundity of $g^x$ with respect to the center $A^x$ (at this
point $A_k$ cannot be replaced by $A_k^+$, since the center $A^x$ may be too
close to the boundary of $A_k^+$ and the rotundity of $g^x$ be very small). 
\subsubsection{}\label{sect2.4.5} We have arrived at the first result of this
paragraph.

Let $S_1,S_2$ be the graphs of $\varphi_j|(0,b],\ j=1,2$ where $\varphi_j$ are
as in section~\ref{sect2.4}, see \eqref{16}. Without loss of generality (see
\cite{HavNer01} or the Introduction) we may assume $A_k\supset S=S_1\cup S_2$.
Suppose $S\subset \{|\zeta|<R\}(=R\D)$.
\begin{theorem}\label{teor6}
For any $k>0$ and any domain $O\supset A_k\cap R\D$, $(S_1,S_2)$ is not a
$bs$-pair in $O$.
\end{theorem}
The proof follows from Lemma~\ref{lema2.2}, since the families $(\phi_j^x)$,
$(K_j^x)$, $j=1,2$, $(g^x)$ enjoying properties (a)--(e) have been constructed
in sections~\ref{sect2.4.1}--\ref{sect2.4.4}.
\subsection{Two graphs in $A_k^+$}\label{sect2.5}
We still assume  $S_1,S_2,\varphi_1,\varphi_2$ are as in section~\ref{sect2.4}
(see \eqref{16}), but we also suppose
\begin{equation}\label{45}
\varphi_1(t)>0,\ t\in(0,b].
\end{equation}
Taking $b<b(k)$ we may assume $S\subset A_k^+$. This time we have to impose a
special proximity conditions on $S_1,S_2$ (\eqref{16} is not sufficient to apply
Lemma~\ref{lema2.2}). Namely, we assume
\begin{equation}\label{46}
\liminf_{x\to 0}\frac {\Delta(x)}{\varphi_1(x)}=0.
\end{equation}
Let $E\subset (0,b]$ be a set with a limit point at the origin and
such that
\[
\lim_{x\to 0,\ x\in E}\frac {\Delta(x)}{\varphi_1(x)}=0.
\]
We define $h(x), X(x)$ (and thus $K_j^x,\ j=1,2$) for $x\in E$ only:
\[
X(x)=x+\frac 12 \varphi_1(x),\ h(x)=2\varphi_1(x)\epsilon(x)
\]
where
\[
\epsilon(x)=\frac{\Delta(x)}{\varphi_1(x)}+\sup\{|\Delta'(t)|:\ 0\le t\le
X(x)\},\ x\in E, 
\]
so that $\lim_{x\to 0,x\in E} h(x)/(X-x)= \lim_{x\to 0,x\in
E}\varepsilon(x)=0$, and we get \eqref{20a}. Now, if $t\in [x,X]$, then, for a
$c\in [x,t]$, $\Delta(t)=\Delta(x)+\Delta'(c)(t-x)\le \epsilon(x)\varphi_1(x)+
\epsilon(x)\varphi_1(x)=h(x)$ which is \eqref{35}. The cell $g^x$ will be
defined for $x\in E$ as the square 
$(x-\varphi_1(x),x+\varphi_1(x))\times (0,2\varphi_1(x))$ with the center
$x+i\varphi_1(x)$. This cell is uniformly rotund.
Moreover, $K_1^x\subset g^x (x\in E)$. Indeed, for $t\in [x,x+\varphi_1(x)]$ we
have $\varphi_1(t)=\varphi_1(x)+\varphi'(c)(t-x)<2\varphi_1(x)$, since we may
assume $|\varphi_1'(c)| <1$. Moreover, $g^x\subset A_k^+$ for any small $x\in
E$, since $\lim_{x\to 0}\frac{2\varphi_1(x)}{x-\varphi_1(x)}=0$.

As in section~\ref{sect2.4.4} we get the following result.
\begin{theorem}\label{teor7}
If $S\subset A_k^+\cap R\D$, and $\varphi_1,\varphi_2$ satisfy \eqref{16} and
\eqref{46}, then, for any domain $O\supset A_k^+;\cap R\D$, $(S_1,S_2)$ is not a
bs-pair in $O$.
\end{theorem}

\section{Concluding remarks}\label{p3}
\subsection{A complete description of some bs-pairs of arcs in
$\C^+$}\label{sect3.1}
Let $\gamma_j:[0,1]\to\C$, $j=1,2$, be two simple
$\mathcal{C}^{1+\varepsilon}$-arcs in $\C$ such that
$\gamma_1(0)=\gamma_2(0)=0$, $\gamma_1(t)\ne \gamma_2(t)$ for $t\in (0,1]$. We
also assume $\Im\gamma_j(t)>0$ for $t\in(0,1]$, $\gamma_j'(t)\ne 0$ for $t\in
[0,1]$.

In this section we define  $S_j$ as the trajectory of $\gamma_j\setminus \{0\}$,
$j=1,2$. Combining Theorem~\ref{teor4}  with
Theorems \ref{teor6} and \ref{teor7} we obtain a complete and very clear
description of all bs-pairs of this sort in $\C^+$. 

Denote by $\tau_j$ the unit tangent vector of $\gamma_j$ at the origin. Note that
$\Im \tau_j\ge 0$, since $S_j\subset \C^+$.
\begin{theorem}\label{teor8}
I. If $\tau_1\ne \tau_2$, then $(S_1,S_2)$ is a bs-pair in $\C^+$. 

II. Suppose $\tau_1=\tau_2$. Then
\begin{itemize}
\item[(a)] $(S_1,S_2)$ is not a bs-pair in $\C^+$, if $\Im\tau_1>0$.
\item[(b)] If $\tau_1$ is real we may assume that $\tau_1=1$, and
for a small $b>0$ the arcs $S_j^b=S_j\cap
\{\Re \zeta<b\}$ for a small $b>0$ are graphs (over $(0,b]$) of real functions
$\varphi_j\in \mathcal{C}^{1+\varepsilon}([0,b])$ as in Theorem~\ref{teor6} (see
\eqref{16}), and $(S_1,S_2)$ is a bs-pair in $\C^+$ if, and only if,
$\liminf_{x\to 0} \frac{\Delta(x)}{\varphi_1(x)}>0$,
$\Delta=|\varphi_1-\varphi_2|$.
\end{itemize}
\end{theorem}
\begin{proof}
In case I the pieces $S_1^b$, $S_2^b$ are separated by two rays in $\C^+$
emanating from the origin, and we may apply Theorem~\ref{teor4} 
of \S\ref{p1} (see also the ``transversal'' example 1 in section~\ref{sect1.4}
and the construction of section~\ref{sect1.6}). In case II(a) $S_1^b$ and
$S_2^b$ stay in an angle in $\C^+$ whose bisector is parallel to $\tau_1$, and
Theorem~\ref{teor6} applies. In case II(b) if \eqref{46} holds, then $(S_1,S_2)$
is not a bs-pair in $\C^+$ by Theorem~\ref{teor7}, since $S_1^b\cup S_2^b$ is
covered by an angle $A_k^+$; if \eqref{46} does not hold, then, taking
$g=(\varphi_1+\varphi_2)/2$ and a small $\mu>0$ (depending on $\limsup_{0}
\Delta/\varphi_1$), $S_1^b$ and $S_2^b$ are separated by the graphs of $g$ and
$(1+\mu)g$, and Theorem~\ref{teor4}  applies; see example 3 in
section~\ref{sect1.4}.
\end{proof}

\subsection{Bounded splittings of functions analytic in the intersection
of Jordan domains}\label{sect3.2}
\begin{figure}
\begin{center}
\psfrag{G1}{$G_1$}
\psfrag{G}{$G$}
\psfrag{N}{$N$}
\psfrag{S}{$S$}
\psfrag{G2}{$G_2$}
\psfrag{1}{$1$}
\psfrag{2}{$2$}
\psfrag{3}{$3$}
\psfrag{4}{$4$}
\psfrag{1p}{$1'$}
\psfrag{2p}{$2'$}
\psfrag{3p}{$3'$}
\psfrag{4p}{$4'$}
\includegraphics{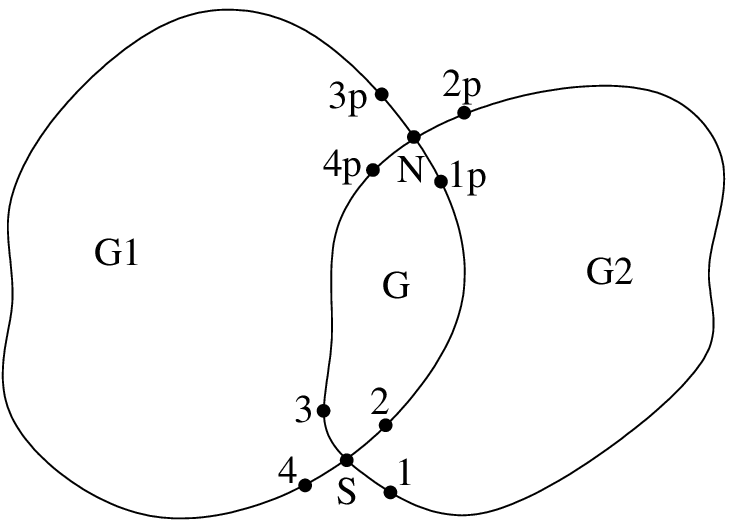}
\end{center}
\caption{}\label{fig3.2}
\end{figure}

Consider two Jordan domains $G_1,G_2$ as on Figure~\ref{fig3.2} and their
intersection $G$.
Suppose the curves $\partial G_1,\partial G_2$ are piecewise
$\mathcal{C}^1$-smooth and intersect \emph{transversally} at $S$ and $N$ (i.e.
any pair of arcs $1,2,3,4$ (respectively $1',2',3',4'$) meet under a positive
angle at $S$, respectively at $N$). It is shown in \cite[Example 4.1, 
section 4.6]{HavNer01} that for any $f\in H^\infty(G)$ there exist functions
$f_j\in H^\infty(G_j)$, $j=1,2$, such that 
\begin{equation}\label{200}
f=f_1+f_2\qquad \text{in } G.
\end{equation}
Now we are going to show that \emph{the transversality assumption cannot be
dropped}.
\begin{figure}
\begin{center}
\psfrag{G1}{$G_1$}
\psfrag{G}{$G$}
\psfrag{A}{$A$}
\psfrag{B}{$B$}
\psfrag{G2}{$G_2$}
\psfrag{ib}{$ib$}
\psfrag{M1}{$M_1$}
\psfrag{M2}{$M_2$}
\psfrag{x}{$x$}
\psfrag{y}{$y$}
\psfrag{C}{$C$}
\psfrag{D}{$D$}
\psfrag{E}{$E$}
\psfrag{F}{$F$}
\psfrag{P}{$P$}
\psfrag{P1}{$P_1$}
\psfrag{Q}{$Q$}
\psfrag{Q1}{$Q_1$}
\psfrag{0}{$o$}
\includegraphics{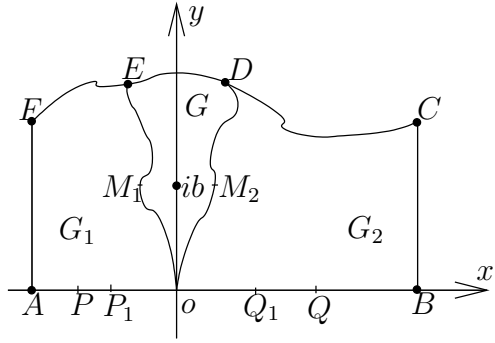}
\end{center}
\caption{$\partial G_1=AoDEFA$, $\partial G_2=oEDCBo$}\label{fig3.2bis}
\end{figure}
Let $G_1,G_2$ be the Jordan domains on Figure~\ref{fig3.2bis}, bounded
respectively, by the loops $AoDEFA$, $oFDCBo$, so that $\partial G=oDEo$.
Suppose $b>0$ is small, and $\partial G\cap\{0\le \Im z \le b\}=
oM_1\cup oM_2$ where the compact arcs $oM_1$, $oM_2$ are as on 
Figure~\ref{fig3.2bis}, so that $(oM_1,oM_2)$ is not a bs-pair with respect to
$\C^+$ (see Theorem~\ref{teor8}). Then \emph{there exists a function $f\in
H^\infty(G)$ which cannot be split as in \eqref{200}}.
\begin{proof}
Take an arbitrary $F\in H^\infty(\C^+\setminus (oM_1\cup oM_2))$ and put
$f=F|G$. Suppose $f$ can be represented by \eqref{200} ($f_j\in H^\infty(G_j),\
j=1,2)$; then we prove
\begin{equation}\label{201}
F=F_1+F_2\quad \text{in }\C^+\setminus (oM_1\cup oM_2)
\end{equation}
with $F_j\in H^\infty(\C^+\setminus oM_j)$, $j=1,2$, a contradiction.

To deduce \eqref{201} from \eqref{200} put $H=G_1\cup G_2$ (so that $\partial
H=AoBCDEFA$) and  note that $f_j$ extends analytically from $G_j$ to $H\setminus
oM_j$, $j=1,2$. Indeed, $f_1=f-f_2$ in $G$ whence $f-f_2$ is the analytic
continuation of $f_1$ from $G$ to $H\setminus oM_1$; the same argument applies
to $f_2$. Extend $f_1,f_2$ to $(H\cup\partial H)'=\C\setminus (H\cup\partial H)$
putting
\[
f_1(\zeta)=f_s(\zeta)=0,\ \zeta\in(H\cup\partial H)',
\]
thus making $f_j\in H^\infty((\partial H\cup oM_j)')$, $j=1,2$. Applying the
preseparation Corollary~3.3 from \cite{HavNer01} we split $f_1,f_2$ in their
domains as follows:
\[
f_j=\varphi_j+r_j,\ j=1,2
\]
where $\varphi_j\in H^\infty([P,Q]\cup o M_j)')$, $r_j\in H^\infty((\partial
H\setminus [P_1,Q_1])')$ (see Figure~\ref{fig3.2bis}). The identity
\begin{equation}\label{202}
F-\varphi_1-\varphi_2=r_1+r_2\quad\text{in }H\setminus (oM_1 \cup oM_2)
\end{equation}
(which is \eqref{200}) shows that $r_1+r_2$ coincides in $H$ with a function
$h\in H^\infty(\C^+)$, since the left side of \eqref{202} is analytic at any
point of $\C^+\setminus H$, and $F,\varphi_1,\varphi_2$ are bounded in their
domains. Therefore \eqref{201} holds with $F_1=\varphi_1,\ F_2=\varphi_2+h$.
\end{proof}
We could get a similar example for the domains $G_1,G_2,G$ as on
Figure~\ref{fig3.2tris}. 
\begin{figure}
\begin{center}
\psfrag{G1}{$G_1$}
\psfrag{G}{$G$}
\psfrag{A}{$A$}
\psfrag{B}{$B$}
\psfrag{G2}{$G_2$}
\psfrag{b}{$b$}
\psfrag{M1}{$M_1$}
\psfrag{M2}{$M_2$}
\psfrag{x}{$x$}
\psfrag{y}{$y$}
\psfrag{C}{$C$}
\psfrag{D}{$D$}
\psfrag{E}{$E$}
\psfrag{F}{$F$}
\psfrag{o}{$o$}
\includegraphics{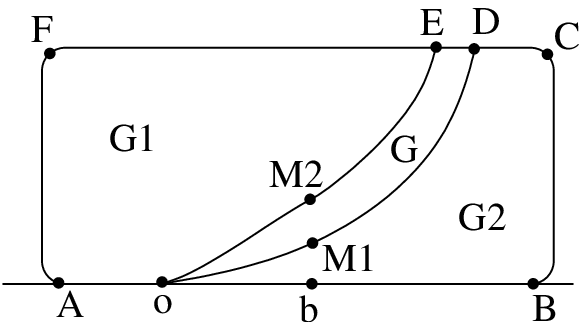}
\end{center}
\caption{$\partial G_1=AoDEFA$, $\partial G_2=oEDCBo$}\label{fig3.2tris}
\end{figure}
The arcs $oM_2$, $oM_2$ are the graphs of functions $\varphi_1,\varphi_2\in
\mathcal C^1([0,b])$ satisfying $\varphi_j(0)=\varphi_j'(0)=0$, $j=1,2$,
$0<\varphi_1(t)<\varphi_2(t)$ ($t\in (0,b]$), $\lim_{t\to 0}
\frac{\varphi_2(t)-\varphi_1(t)}{\varphi_1(t)}=0$ (see Theorem~\ref{teor8}).
However this kind of tangency (when $oM_1$, $oM_2$ are tangent to $\partial H$
at the origin) is compatible with the splitting formula \eqref{200} for
\emph{any} $f\in F^\infty(G)$ if $\lim_{t\to 0}\frac{\varphi_2(t)-\varphi_1(t)}
{\varphi_1(t)}>0$. This happens, for example, for 
\[
G_1=\{0<x<1,\ 0<y<2x^2\},\qquad G_2=\{0<x<1,\ x^2<y<2\}.
\]

\subsection{bs-pairs of hyperbolically close sets}\label{sect3.3}
Very satisfactory looking Theorem~\ref{teor8} deals with pairs of \emph{graphs}
and cannot be applied to disconnected sets. In this section we describe examples
of bs-pairs $(S_+,S_-)$ with respect to the \emph{right} half plane $\Pi$ with
$S_+$ and $S_-$ hyperbolically very close.

Let $g$ be a non-negative function on $[0,b]$ such that $g(\xi)\le \xi$
($\xi\in[0,b]$) and $g(\xi)>0$ for $\xi\in(0,b]$. Consider a \emph{strictly}
decreasing sequence $(\xi_n)_{n=1}^\infty$ in $(0,b]$ tending to zero and put
$\zeta_n=\xi_n+ig(\xi_n)$. For $r>0$ denote by $B_n(r)$ the \emph{closed} disc
$\{|\zeta-\zeta_n|\le r\}$. Suppose a sequence $(r_n)_{n=1}^\infty$ of positive
numbers satisfies
\begin{equation}\label{100}
\sum_{n=1}^\infty \frac {r_n}{g(\xi_n)}<\infty,
\end{equation}
so that $\sum r_n/\xi_n <\infty$ as well. We assume $r_n<g(\xi_n)$ for all $n$
and $r_n<\min(\xi_n-\xi_{n-1},\xi_{n+1}-\xi_{n})$, so that the discs
$B_n=B_n(r_n)$ are disjoint and stay in $\C^+$. Put
\[
S_+= \cup_{n=1}^\infty B_n,\quad S_-= \cup_{n=1}^\infty \overline{B_n},\qquad 
S=S_+\cup S_-.
\]

\begin{figure}
\begin{center}
\psfrag{z1}{$\xi_1$}
\psfrag{z2}{$\xi_2$}
\psfrag{z3}{$\xi_3$}
\psfrag{zn}{$\xi_n$}
\psfrag{g(zn)}{$g(\xi_n)$}
\includegraphics{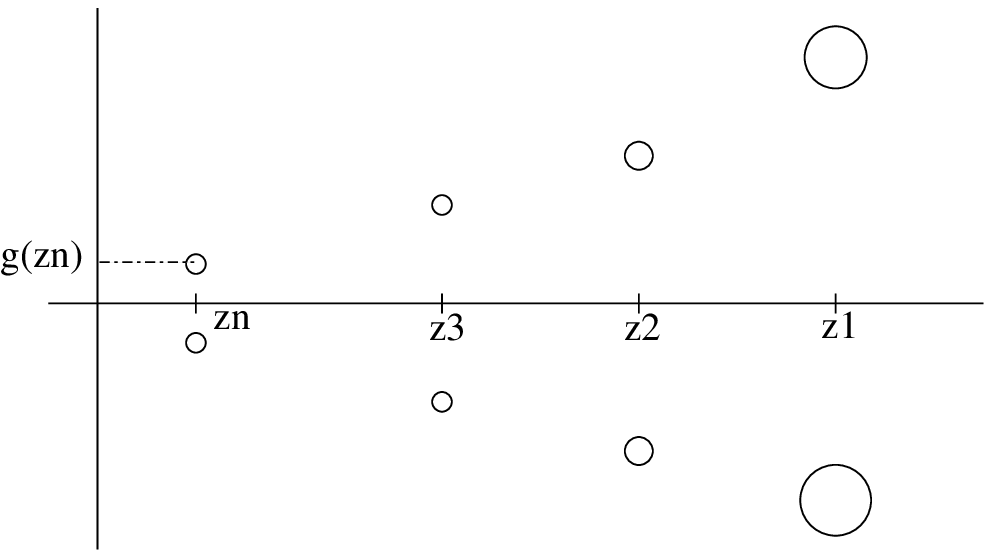}
\end{center}
\caption{}\label{fig3.3}
\end{figure}

\begin{theorem}\label{teor9}
$(S_+,S_-)$ is a bs-pair with respect to the right half plane $\Pi$.
\end{theorem}
Note that the hyperbolic distance between $\xi_n$ and $\bar \xi_n$ (and $B_n$
and $\overline{B_n}$) in $\Pi$ is comparable with $g(\xi_n)/\xi_n$ and tends to
zero, if say, $g'(0)=0$.
\begin{proof}
Any $f\in H^\infty(\Pi\setminus S)$ has angular boundary values a.e. on $i\R$
and $\partial S$, and can be represented as follows:
\[
f=f_1+f_2=f_1+f_++f_-\qquad \text{in }\Pi
\]
where 
\begin{equation}\label{101}
f_1=-(1+z)C_{i\R}^{f/(1+z)},\qquad
f_+=(1+z)\sum_{n=1}^\infty C_{\partial B_n}^{f/(1+z)},\qquad
f_-=(1+z)\sum_{n=1}^\infty C_{\partial \overline{B_n}}^{f/(1+z)}.
\end{equation}
($C_A^k$ is the Cauchy type integral $\frac 1{2\pi i}\int_A
\frac{k(z)}{z-\zeta}dz$, $\zeta\notin A$).
The imaginary axis is oriented ``upwards'', and the circles $\partial B_n$,
$\partial \overline{B_n}$ are oriented clockwise. The proof of \eqref{101} is
standard ($f(iy)/(1+y^2)\in L^1(\R,dy)$, $f(z)/(1+z)=O(1/z),\ |z|\to\infty$).
Now, $f_\pm = (1+z)C^{\mu_\pm}$ where $\mu_{\pm}$ are \emph{finite} complex
charges on $S_\pm$, so that $f_2$ is defined and analytic in $\C\setminus(S\cup
\{0\})$ and $f_2(z)=O(1)$ ($|z|\to\infty$). Clearly $f_1=f-f_2$ is analytic in
$\Pi$, and $f_1(z)=O(1)$ ($|z|\to\infty$).

We are going to prove that 
\begin{equation}\label{102}
f_1\in H^\infty(\Pi)
\end{equation}
whence $f_2\in H^\infty(\Pi\setminus S)$. At last we prove $f_\pm$ are bounded
in $\Pi\setminus S_\pm$ and get the final splitting of $f$ into the sum of
elements of $H^\infty(\Pi\setminus S_\pm)$:
\[
f=(f_1+f_+)-f_-.
\]
(a) $f_+$ is bounded on any set $\C\setminus A_k^+$, $k>2$, where 
$A_k^+=\{\xi+i\eta:\ \xi>0, 0<\eta< k\xi\}$ (according to our assumptions
$S_+\subset A_2$). Indeed,
\begin{equation}\label{104}
|C_{\partial B_n}^f(\zeta)|\le \|f\|_\infty r_n / \dist(\zeta,B_n),\qquad
\zeta\notin S_+.
\end{equation}
Let $\zeta\in A_k^+$. If $\Re \zeta\le 0$, then $\dist(\zeta,B_n)\ge
\xi_n-r_n\ge g(\xi_n)-r_n$; if $\Im \zeta\le 0$, then $\dist(\zeta,B_n)\ge
g(\xi_n)-r_n$, and if $\zeta\in \Pi\setminus A_k^+$, $\Im\zeta>0$, then 
\[
\dist(\zeta,B_n)\ge |\zeta-\zeta_n|-r_n\ge |\zeta_n|\sin\varphi-r_n\ge
g(\xi_n)\sin\varphi -r_n,
\]
where $\varphi$ is the angle in $\Pi$ with sides $\eta=2\xi,\ \eta=k\xi$. Thus,
for a $c_k>0$, 
\begin{equation}\label{105}
\dist(\zeta,B_n)\ge c_k g(\xi_n)
\end{equation}
for all $\zeta\in \C^+\setminus A_k^+$ and $n=1,2,\ldots$. Combining
\eqref{104}, \eqref{105} and \eqref{100} we see that $f_+$ is analytic
and bounded in
$\C\setminus A_k^+$, $k>2$. The same proof applies to $f_-$, so that $f_-$ is
bounded in $\C\setminus \overline{A_k^+}$, $k>2$. Putting $A_k= A_k^+\cup 
\overline{A_k^+}\cup (0,+\infty)$ we conclude that $f_2$ is bounded in
$\C\setminus A_k$, and $f_1=f-f_2$ \emph{is bounded in} $\Pi\setminus A_k$. But
$f_1(\zeta)=O(1)$ ($|\zeta|\to\infty,\zeta\in \Pi$), and we may fix a large
$R>0$ making $f_1$ bounded in $\Pi\setminus s$ where $s=A_k\cap \{|\zeta|<R\}$
(we fix a $k>2$, say $k=3$). It remains to estimate $f_1$ in $s$. Being analytic
in $s$, $f_1$ is bounded on its sides and the arc $A_k\cap\{|\zeta|=R\}$ whereas
\begin{equation}\label{106}
|f_1(\zeta)|=O(\log|\zeta|),\ |\zeta|\to 0, \zeta\in s
\end{equation}
and (a very weak form of) Phragmen-Lindel\"of applies whence $f_1$ \emph{is
bounded in} $s$, and \eqref{102} holds. (To get \eqref{106} write (using
$|\eta|<k\xi$)
\[
\begin{split}
|f_1(\zeta)|\lesssim & \int_{-\infty}^{+\infty} \frac{dy}{(1+|y|)(\xi+|\eta-y|)}
\lesssim \int_{|y|<2k\xi}+ \int_{|y|>2k\xi}=\\
=& O(1)+ O\left(\int_{2k\xi}^{+\infty}
\frac{dy}{(1+y)y}\right)= O(|\log\xi|)=O(|\log\zeta|),
\end{split}
\]
since $|\zeta|\simeq \xi$ for $\zeta\in s$).

(b) From (a) we see that $f_2=f_+ + f_-$ is bounded in $\Pi\setminus S$ (and, in
fact, in $\C\setminus(S\cup\{0\})$). We already know $f_+$ is bounded in
$\C\setminus (A^+_k\cup\{0\})$. We only have to show $f_+$ is bounded in
$A_k^+\setminus S_+$. But in $A_k^+\setminus S_+$
\[
f_+=f-f_1-f_{-};
\]
$f$ is bounded in $\Pi\setminus S$, $f_1$ is bounded in $\Pi$ (see \eqref{102})
and we have proved $f_-$ is bounded in $\C^+$. Thus $f_+\in
H^\infty(A_k^+\setminus S_+)$, $f_+\in H^\infty(\C\setminus(S^+ \cup\{0\}))$.
The boundedness of $f_-$ in $\Pi\setminus S_-$ follows now from $f_-=f_2-f_+$.
\end{proof}

\def\cprime{$'$}
\providecommand{\bysame}{\leavevmode\hbox to3em{\hrulefill}\thinspace}
\providecommand{\MR}{\relax\ifhmode\unskip\space\fi MR }
\providecommand{\MRhref}[2]{%
  \href{http://www.ams.org/mathscinet-getitem?mr=#1}{#2}
}
\providecommand{\href}[2]{#2}

\end{document}